\documentclass[doublespacing]{article}
\usepackage[utf8]{inputenc}
\usepackage[left=4.00cm, right=3.50cm, top=2.50cm, bottom=2.50cm]{geometry}
\usepackage{amsmath}
\usepackage{amsfonts}
\usepackage{amssymb}
\usepackage{graphicx}
\usepackage{setspace}
\usepackage{colortbl}
\usepackage{wrapfig}
\usepackage[thmmarks,amsmath]{ntheorem}
\usepackage{tikz-cd}
\usepackage[colorlinks,allcolors=blue]{hyperref}

\theoremheaderfont{\bfseries\upshape}
\theorembodyfont{\normalfont}
\theoremstyle{plain}

\newtheorem{definition}{Definition}[section]
\newtheorem{theorem}{Theorem}[section]
\newtheorem{remark}{Remark}[section]
\newtheorem{proposition}{Proposition}[section]
\newtheorem{corollary}{Corollary}[section]
\newtheorem{example}{Example}[section]
\newtheorem{examples}{Examples}[section]
\newtheorem{lemma}{Lemma}[section]
\newcommand{\proof}{\begin{bfseries} \fontsize{3mm}{3mm} \selectfont Proof \end{bfseries} \\}

\newcommand{\vs}{\vspace{5mm}}

\author{Federica Pasqualone}
\title{Gabriel-Ulmer Duality for Topoi, An Introduction}
\date{June 5, 2024}

\begin{document}

\maketitle
\vs 

\abstract{The aim of this paper is to relate the classical result of Gabriel-Ulmer to the geometry of topoi. The usage of the attribute 'left exact' when dealing with functors involved in this duality is indeed not casual and it is related to the geometrical side of the story, i.e. sheaf theory, as in \cite{SGA}. Thought to be a very basic introduction to the subject, it is mostly self-contained. The reader is assumed to be familiar with the fundamentals of category theory, no further prerequisite knowledge is required. The proof of Gabriel-Ulmer duality follows \cite{MP}.}

\vs 
\tableofcontents

\newpage 

\section{Introduction}\label{intro}

\ \ \ Gabriel-Ulmer duality is a well-established result in categorical logic stating the equivalence between essentially algebraic theories and locally finitely presentable categories. Essentially algebraic theories extend the class of Lawvere algebraic theories by allowing for partially defined operations. In categorical terms, they are represented by small categories having all finite limits, i.e. finite products and equalizers.\\

\begin{example}
   The theory of categories is essentially algebraic. In fact, composition is only allowed when the target of the first arrow corresponds to the source of the second.
\begin{center}
   \begin{tikzcd}
    A \arrow[r,"f"] &B \arrow[r, "g"] &C
\end{tikzcd}
\end{center}
\vs 
\end{example} 

The project of this paper started after my talk ''The Dark Side of Gabriel-Ulmer Duality'', where I explained why we call the relevant functors involved in the mathematical formalization of this duality 'left exact'. Indeed, the link to homological algebra, or better to sheaf theory, is apparent in \cite{SGA}, but it is somehow made obscure in the more logic-based textbooks. An extensive answer is provided in \ref{appl}. \\

This paper develops the basics, giving large space to preliminary results and tools, a more advanced account of the relation with topos theory is provided by Ivan Di Liberti and Julia Ramos González in \cite{DLRG}. Moreover, Gabriel-Ulmer duality can be further extended to include clans, small categories equipped with a particular class of maps called display maps, as in the work of Jonas Frey \cite{JF}.\\

Future work \ref{future} will develop the geometry of clans. The claim is that clans generalize the notion of small sites, where display maps play the role of coverings.

\section{Preliminaries}\label{pre}

\ \ \ Let us recall facts and definitions that are crucial for the subject starting from the concept of Lawvere algebraic theory. 

\subsection{Algebraic Theories}\label{AT}

\begin{definition} An algebraic theory is a small category with finite products and it is usually denoted by $\mathcal{T}$.
\end{definition}

\begin{example}
The theory of Abelian groups is algebraic. Objects are natural numbers, morphisms $\phi: m \rightarrow n$ between them are represented by $n \times m$ matrices.
\end{example}
\begin{example}
For $R$ a unital ring,  $\text{Mod}_{R}$, the category of $R$-modules, is small and with finite products. The theory of (right) $R$-modules is therefore algebraic.
\end{example}
\vs 
\begin{definition} An algebraic theory $\mathcal{T}$ is called essentially algebraic if the operation of composition is partially defined.
\end{definition}
\begin{example} Horn theories are essentially algebraic, the viceversa is not necessarily true. An Horn theory in a many-sorted-first-order signature $\Sigma$ is a set of Horn sequents $\phi \vdash \psi$ in a certain context $\vec{x}$, possibly empty. The formulas appearing in an Horn sequent are Horn formulas, meaning they are derived from atomic formulas and $\top$ using only $\land$, no other conjunction is allowed.
For more, \cite{PV}.
\end{example}

\subsection{Left-exactness}\label{LE}

\ \ \ Since we will deal mostly with two categories, recall that a 2-category \cite{KE} has categories as objects, functors between them as morphisms - defining a 1-cell- and natural transformation between such functors - creating a two cell. Cells can be further composed both horizontally and vertically and the two kinds of composition satisfy some coherence conditions. 
\vs
\begin{center}
    \begin{tikzcd}
        \mathbb{C} \arrow[bend left]{r}[name=F]{F} \arrow[bend right ]{r}[name=G, swap]{G} & \mathbb{D}
        \arrow[phantom, from=F, to=G, "\Rightarrow" rotate= -90]
        \end{tikzcd}
\end{center}
\vs

  Relevant 2-categories:
\begin{itemize}
    \item[(i)]\textbf{Lex}, the 2-category of small categories with all finite limits and finite limits preserving functors;
    \item [(ii)]\textbf{LEX}, the 2-category of categories with all finite limits and finite limits preserving functors;
    \item [(iii)] \textbf{LFC}, the 2-category of categories having all (small) limits and filtered colimits with limits and filtered colimits preserving functors;
    \item [(iv)] \textbf{LFP}, the 2-category of locally finitely presentable categories with limits and filtered colimits preserving functors;
    \item[(v)] \textbf{CAT}, the 2-category of categories, with functors as 1-cells and natural transformations between them as 2-cells.
\end{itemize}

\ \ \ A model of an essentially algebraic theory $\mathbb{C} \in \mathbf{Lex}$ is given by a lex (left exact) functor $ F \in \mathbf{LEX}\left(\mathbb{C}, \mathbf{Set}\right)$. \\

The term left exact functor sounds familiar to the ones with some knowledge of homological algebra, and indeed the definition
\begin{definition}{[\cite{SGA},1, \S2]} Let $\mathbb{C}$ be a category with finite projective limits, a functor $F: \mathbb{C} \rightarrow \mathbb{C}'$ is left exact iff it commutes with finite projective limits.
\end{definition}

\noindent extends the same notion given for an additive functor of Abelian categories 
\begin{definition}
    A covariant additive functor $\Phi: \mathbb{A} \rightarrow \mathbb{A}'$ is  left exact iff whenever $0 \rightarrow A \rightarrow B \rightarrow C$ is exact, $0 \rightarrow \Phi A \rightarrow \Phi B \rightarrow \Phi C$ is exact.
\end{definition}
\begin{remark} An Abelian category is a category enriched over Abelian groups. It has a zero object (an object that is both initial and terminal), all binary products, kernels and cokernels and all the monics and epics are normal (kernels and co-kernels, respectively, of some morphisms).
\end{remark}

The link with homological algebra for the left exactness property just displayed will become evident in discussing presheaves (and sheaves) in the section on topos theory - and in \cite{SGA}. 
\begin{remark}
    The category of left exact functors with values in $\mathbf{Set}$ provides also models for lim-theories, first-order theories over a language $\mathcal{L}$ with axioms of the form 
    \begin{equation*} 
        \forall x \left( \phi \left(x \right) \rightarrow \exists! y \psi \left(x,y\right)\right)
    \end{equation*} where the formulae $\phi$ and $\psi$ are finite conjunctions of atomic ones.
\end{remark}

\subsection{Filtered categories and colimits}

\begin{definition}{[\cite{CWM}, IX, \S1]} A category $\mathbb{J}$ is said to be filtered iff it is non empty and the following applies:
    \begin{itemize}
        \item[(i)] For any two objects $i, j \in \mathbb{J}_0$, there exists $ k \in \mathbb{J}_0$ and morphisms $f_i,f_j \in \mathbb{J}_1$ as displayed below
        \vs 
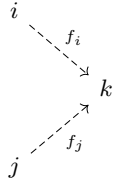
\begin{center}
        \begin{tikzcd}
            i \arrow[dr, dashed, "f_i"] &{}\\
            {} &k \\
            j \arrow[ur, dashed,swap, "f_j"] &{}\\
        \end{tikzcd} 
    \end{center}
    \vs
    \item[(ii)] For any two parallel pair of morphisms \begin{tikzcd}[column sep=.5cm] f,g: h \arrow[r, shift right] \arrow[r, shift left] & j \end{tikzcd}, there exists an object $k \in \mathbb{J}_0$ and an arrow $w: j \rightarrow k$ such that the diamond commutes
    \vs 
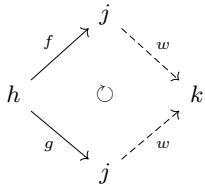
    \begin{center}
        \begin{tikzcd}
            {} &j \arrow[dd, phantom, "\circlearrowright"] \arrow[dr , "w", dashed] &{}\\
            h \arrow[ur, "f"] \arrow[dr, swap, "g"] &{} &k\\
            {} &j \arrow[ur, "w", swap, dashed] &{}
        \end{tikzcd}
    \end{center}
    \vs

    \end{itemize}
\end{definition}
\begin{remark}
    The notion of filtered, or co-filtered, category - depending on the authors - generalizes the one of directed preorder, i.e. a preorder in which every element has an upper bound.
\end{remark}
\begin{definition}{[ibid.]} A filtered colimit is a colimit of a diagram $D: \mathbb{J} \rightarrow \mathbb{C}$ over a filtered category.
\end{definition}  

\begin{remark} The category of sets, $\mathbf{Set}$, is complete and co-complete and finite limits commute with limits and filtered colimits. These properties make it particularly suitable for the development of the theory.
\end{remark}

\ \ \ The category $\mathbf{LEX}\left(\mathbb{C}, \mathbf{Set}\right)$ has all (small) limits and filtered colimits and they are preserved by the evaluation functor
\begin{align*}
    \mathbf{LEX}\left(\mathbb{C}, \mathbf{Set}\right) &\rightarrow \mathbf{Set}\\
    F &\mapsto F\left(C\right)
\end{align*}
for $C \in \mathbb{C}_0$.

\subsection{Locally Finitely Presentable Categories}\label{LFPC}

\ \ \ Varieties of algebras, domains, lattices, are some instances of locally presentable categories. It is a quite comprehensive class, however there is even a more general account provided by accessible categories that allows to describe fields, Hilbert spaces and other interesting mathematical structures. \\

Gabriel-Ulmer duality investigates in particular locally \textit{finitely} presentable categories, categories with all colimits, whose objects can be obtained as directed colimits of objects in a set. They generalize complete algebraic lattices.
\begin{definition} A diagram over a category $\mathbb{C}$ is said to be directed iff it is indexed by a directed poset $\left(I, \le\right)$  - thought of as posetal category.
\end{definition}
\begin{remark} A poset $\left(I, \le\right)$ is said to be directed iff each pair of elements $x, y \in \left(I, \le \right)$ has an upper bound. To make a category out of a poset we define arrows as follows:
    \begin{equation*}
        x \le y \iff x \rightarrow y
    \end{equation*}
    Therefore, each hom-set of a so-called thin (or posetal) category has at most one arrow.
\end{remark}
\begin{definition} Colimits of a directed diagram are called directed colimits or, equivalently, left limits, projective limits.
\end{definition}

\begin{definition}{[\cite{AR},1,1.A, Definition 1.1]} An object $A \in \mathcal{A}_0$ of a locally small category is called finitely presentable (f.p.) iff the functor
    \begin{equation*}
        \mathcal{A}\left(A, -\right): \mathcal{A}\rightarrow \mathbf{Set}
    \end{equation*}
    preserves filtered colimits. The full subcategory of finitely presentable objects of $\mathcal{A}$ is denoted by $\mathcal{A}_{f.p}$.
\end{definition}

\begin{examples}
    \begin{itemize}
        \item[(i)] $K \in \mathbf{Set}_0$ is finitely presentable iff $K$ is finite;
        \item[(ii)] Let S be a set of sorts, then an S-sorted set is representable iff it has finite power. Given an object $X \in \mathbf{Set^S}, \ X = \left\{X_s\right\}_{s \in S}$, its power is defined as 
    \begin{equation*} |X| = \sum_{s \in S} |X_s| \end{equation*}
        \item[(iii)] A group $G \in \mathbf{Grp}_0$ is finitely presentable iff can be represented by a finite number of generators and relations, i.e. as quotient of the free group on n-generators by the relations. For instance, the additive group $\mathbb{Z}$ of integers is finitely presentable, $\mathbb{Z}_n$ is also finitely presentable (being obtained from the previous one by adding the congruence relation $x = nz + y$, for some $z \in \mathbb{Z}$). Therefore, every finite cyclic group is. 
        \item[(iv)] Let $\mathbb{C}$ be small, every hom-functor is finitely presentable in $\mathbf{Set}^{\mathbb{C}}$. Indeed, the Yoneda embedding
        \begin{align*}
            y :\mathbb{C} &\rightarrow \mathbf{Set}^{\mathbb{C}^{\text{op}}}\\
            K &\mapsto \mathbb{C}\left(-, K \right)\\
            \mathbb{C}\left(K,K'\right) \ni f  &\mapsto \mathbb{C}\left(-,f\right):\mathbb{C}\left(-, K\right) \rightarrow \mathbb{C}\left(-, K'\right)
        \end{align*}
        is representable.
        \item[(v)]  A topological space is finitely presentable in $X \in \mathbf{Top}_0$ iff X is finite and discrete. 
    \end{itemize}
\end{examples}

\begin{proposition}{[\cite{AR},1, 1.A, Proposition 1.3 ]} A finite colimit of finitely presentable objects is finitely presentable.
\end{proposition}

\begin{definition}{[\cite{AR}, 1, 1.A, Definition 1.9]} A category $\mathcal{A}$ is called locally finitely presentable (l.f.p.) iff it is co-complete and it has a collection $\mathcal{G}$ of finitely presentable objects: every object is a projective limit of objects of $\mathcal{G}$. The elements of such collection are called generators of $\mathcal{A}$.
\end{definition}
By definition, $\mathbf{LFP}\subseteq \mathbf{LFC}$. The locally finitely presentable categories are a full subcategory of the ones having all finite limits and filtered colimits. In fact, the previous definition can be restated as follows:

\begin{definition}
A category $\mathcal{A}$ is called locally finitely presentable iff $\mathcal{A} \in \textbf{LFC}_0, \mathcal{A} $ is locally small and it has a set of generators consising of finitely many presentable objects.
\end{definition}
\begin{examples}
    \begin{itemize}
        \item[(i)] $\mathbf{Set}$ is locally finitely presentable, since each set is a directed colimit of the diagram of all its finite subsets ordered by inclusion and there exists an essentially unique countable set of finite sets;
        \item[(ii)] A thin category is locally finitely presentable iff it is a complete algebraic lattice;
        \item[(iii)] Non-examples: $\mathbf{Set}^\text{fin}$ is not, indeed is not co-complete, $\mathbf{Top}$.
    \end{itemize}
\end{examples}

\begin{remark}
    The different notations for categories, $\mathbb{C}$ and $\mathcal{A}$, is meant to facilitate the distinction between objects of $\mathbf{LEX}$ and objects of $\mathbf{LFC}$.
\end{remark}

\subsection{Towards the Proof of Gabriel-Ulmer Duality}\label{tow}

\ \ \ The following inclusions hold:
\begin{align*}
    \mathbf{LEX} &\hookrightarrow \mathbf{CAT}\\
    \mathbf{LFC} &\hookrightarrow \mathbf{CAT}
\end{align*}
these functors are faithful, full on 2-cells.\\

\begin{remark}
    The bar for unit and co-unit displayed below is a useful trick for remembering the order of composition of the functors to point towards the identity over the right category.
\end{remark}

The category of $\mathbf{Set}$ not only is an object of the both subcategories, $\mathbf{Set} \in \mathbf{LEX}_0 \cap \mathbf{LFC}_0$, but the $\mathbf{LEX}$ operations (finite limits) commute with the $\mathbf{LFC}$ ones (small limits and filtered colimits). In other words, [MP] $\mathbf{Set}$ is a symmetric $\left(\mathbf{LEX}, \mathbf{LFC}\right)$- bistructure, namely there exists a 2-adjoint pair of 2-functors
\begin{center}
    \begin{tikzcd}
    \epsilon|\mathbf{LEX}^{\text{op}} \arrow[r, shift left, "G"] &\mathbf{LFC} \arrow[l, shift left,"F"] |\eta 
\end{tikzcd}
\end{center}
This is the crucial adjunction for the proof of the Gabriel-Ulmer duality theorem. Since the 2-categories of interest are subcategories of $\mathbf{CAT}$, it is reasonable to construct the adjunction on  $\textbf{CAT}$ and then study its restriction.\\
On the 2-category of categories we have
\begin{center}
    \begin{tikzcd}
    \epsilon_0| \mathbf{CAT}^{\text{op}} \arrow[r, shift left, "G_0"] &\mathbf{CAT} \arrow[l, shift left,"F_0"]| \eta_0
\end{tikzcd}
\end{center}
where
\begin{align*}
    G_0= \textbf{CAT}\left(-, \mathbf{Set}\right) \ \ \ && \left(\epsilon_0\right)_{\mathbb{C}}: \mathbb{C} \rightarrow \mathbf{CAT}\left(\mathbf{CAT}\left(\mathbb{C}, \mathbf{Set}\right), \mathbf{Set}\right)\\
    F_0= \textbf{CAT}\left(-, \mathbf{Set}\right) \ \ \ && \left(\eta_0\right)_{\mathcal{A}}: \mathcal{A} \rightarrow \mathbf{CAT}\left(\mathbf{CAT}\left(\mathcal{A}, \mathbf{Set}\right), \mathbf{Set}\right)
\end{align*}
i.e. the components of unit and co-unit are just given by evaluations.\\

$\mathbf{Set}$ is symmetric, therefore:
\begin{itemize}
    \item[(i)] $\forall \mathbb{C} \in \mathbf{LEX}_0, \ \mathbf{LEX}\left(\mathbb{C}, \mathbf{Set}\right) \in \mathbf{LFC}$;
    \item[(ii)] $\forall \mathcal{A} \in \mathbf{LFC}_0, \ \mathbf{LFC}\left( \mathcal{A}, \mathbf{Set}\right) \in \mathbf{LEX}$.
\end{itemize}
 
Any morphism $\phi \in \mathbf{LEX}\left(\mathbb{C}, \mathbb{C}'\right)$ induces a morphism $\hat{\phi} \in \mathbf{LFC}\left( \mathbf{LEX} \left(\mathbb{C}', \mathbf{Set}\right), \mathbf{LEX} \left(\mathbb{C},\mathbf{Set}\right) \right)$, and dually for $\mathbf{LFC}_1.$
This implies the existence of 2-functors
\begin{align*}
            G: \mathbf{LEX}^{\text{op}} &\rightarrow \mathbf{LFC} \hfill & F:\mathbf{LFC} &\rightarrow \mathbf{LEX}^{\text{op}}\\
            \mathbb{C} &\mapsto \mathbf{LEX}\left(\mathbb{C}, \mathbf{Set}\right) \hfill  & \mathcal{A} &\mapsto \mathbf{LFC}\left(\mathcal{A}, \mathbf{Set}\right)
        \end{align*}

        We still need to define unit and co-unit, namely:
        \begin{align}
            &\epsilon: FG \rightarrow \text{id}_{\mathbf{LEX}^{\text{op}}}\\
            &\eta: \text{id}_{\mathbf{LFC}^{\text{op}}} \rightarrow GF
        \end{align}

        For (1), since $\mathbf{LEX} \hookrightarrow \mathbf{CAT}$,
        \begin{equation*}
            \mathbf{CAT}\left(\mathbf{CAT}\left(\mathbb{C}, \mathbf{Set}\right), \mathbf{Set}\right) \rightarrow \mathbf{CAT}\left(\mathbf{LEX}\left(\mathbb{C}, \mathbf{Set}\right), \mathbf{Set}\right)
        \end{equation*}
        and, precomposing with $\left(\epsilon_0\right)_{\mathbb{C}}$, the resulting functor factors through
        \begin{equation*}
            \mathbf{LFC}\left(\mathbf{LEX}\left(\mathbb{C}, \mathbf{Set}\right), \mathbf{Set}\right) \rightarrow \mathbf{CAT}\left(\mathbf{LEX}\left(\mathbb{C}, \mathbf{Set}\right), \mathbf{Set}\right)
        \end{equation*}
        yielding
        \begin{equation*}
            \epsilon_{\mathbb{C}}: \mathbb{C} \rightarrow \textbf{LFC}\left(\mathbf{LEX}\left(\mathbb{C}, \mathbf{Set}\right), \mathbf{Set}\right)
        \end{equation*}
        Similarly, for the components of $\eta.$\\

        Further, we observe that 
        \begin{align}
            \mathbf{Lex} &\hookrightarrow \mathbf{LEX}\\
            \mathbf{LFP} &\hookrightarrow \mathbf{LFC}
        \end{align}
        as full sub-2-categories.  Equation (3) is valid because $\mathbf{Lex}$ is the category of essentially small objects of $\mathbf{LEX}$.
        To better formalize it, we introduce the following
        \begin{definition} A sub-category $\mathbb{D}' $ of a category $\mathbb{D}$ is called closed under isomorphisms, or replete, iff
            for every isomorphism $f\in \mathbb{D}\left(D_1,D_2\right)$ such that $D_1 \in \mathbb{D}'_0$, $f \in \mathbb{D}_1'$.\\ In other words, it is replete iff  
            every isomorphism of the mother category that has source in the sub-category is an arrow in the sub-category.
            Thus, also the target and the inverse of $f$ belong to the sub-category.
        \end{definition}
    \ \ \ Back to the case of interest, an essentially small object is isomorphic to a small category. Therefore, for any essentially small object $\mathbb{D} \in \mathbf{LEX}_0$ there exists an isomorphism 
    \begin{equation*}
        \mathbb{C} \overset{h}{\rightarrow} \mathbb{D}
    \end{equation*}
    where $\mathbb{C}$ is a small category, i.e. $\mathbb{C} \in \mathbf{Lex}$, replete sub-category of $\mathbf{LEX}$.

\subsection{Proof of the Theorem}\label{pr}

    \ \ \ We are now in position to give a proof of the Gabriel-Ulmer Duality Theorem:
    \begin{theorem}{[\cite{MP}, Theorem 1.2]} The adjunction  $\begin{tikzcd}
        \epsilon|\mathbf{LEX}^{\text{op}} \arrow[r, shift left, "G"] &\mathbf{LFC} \arrow[l, shift left,"F"] |\eta 
    \end{tikzcd}$ restricts to an equivalence of 2-categories
    \begin{center}
        \begin{tikzcd}
            \epsilon|\mathbf{Lex}^{\text{op}} \arrow[r, shift left, "G"] &\mathbf{LFP} \arrow[l, shift left,"F"] |\eta 
        \end{tikzcd}
    \end{center}
Namely,
\begin{itemize}
    \item[(i)] If $\mathbb{C} \in \mathbf{Lex}_0$, $G\left(\mathbb{C}\right) \in \mathbf{LFP}_0$ and $\epsilon_{\mathbb{C}}$ is an equivalence of categories;
    \item[(ii)] If $\mathcal{A} \in \textbf{LFP}_0, \ F\left(\mathcal{A}\right) \in \mathbf{Lex}_0$ and $\eta_{\mathcal{A}}$ is an equivalence.
\end{itemize}
\end{theorem}

\begin{lemma}{[\cite{MP}, Lemma 1.4]} Let $\mathcal{A}$ be a locally small category with all limits and filtered colimits preserved by the functor $G: \mathcal{A} \rightarrow \mathcal{B}$, for $\mathcal{B} \in \mathbf{LFP}_0.$ If $\mathcal{A}$ has a set of generators $\mathcal{G}_{\mathcal{A}}$, there exists a functor $F:\mathcal{B}\rightarrow \mathcal{A}$ such that $G \vdash F$.
\end{lemma} 

\begin{corollary}{[\cite{MP}, Corollary 1.5]} Let $\mathcal{A} \in \mathbf{LFP}, A \in \mathcal{A}_{f.p}$, then $\mathcal{A}\left(A, - \right) \in \mathbf{LFC}\left(\mathcal{A}, \mathbf{Set}\right).$ Therefore, there exists a canonical functor
    \begin{align*}
        h_{\mathcal{A}}: \mathcal{A}_{f.p}^{\text{op}} &\rightarrow \textbf{LFC}\left(\mathcal{A}, \mathbf{Set}\right)\\
A &\mapsto \mathcal{A}\left(A, - \right)
    \end{align*}
    that is an equivalence of categories.
\end{corollary}

\proof 
     By the Yoneda Lemma, $h_{\mathcal{A}}$ is fully faithful. To show: $h_{\mathcal{A}}$ essentially surjective.\\  Consider a functor $X \in \textbf{LFC}\left(\mathcal{A}, \mathbf{Set}\right)$, from the previous lemma it follows that it has a left adjoint $X \vdash Y$. Therefore,
    \begin{equation*}
        \mathcal{A}\left(Y\left(1\right), A \right) \cong \mathbf{Set}\left(1, XA\right) \cong XA \ \ \ \forall A \in \mathcal{A}_0
    \end{equation*} 
    Thus, there exists a natural iso $\mathcal{A}\left(Y\left(1\right), - \right) \cong X$. Since X preserves filtered colimits $Y\left(1\right)$ is finitely presentable by definition. $\square$ 
    
    \proof  
    \textbf{Theorem 1.1 (i)} Let $\mathcal{A} \in \textbf{LEX}\left(\mathbb{C}, \mathbf{Set}\right)$, $\mathcal{A}$ is l.f.p.., therefore $\mathcal{A}$ is representable. \\
    The Yoneda embedding $y: \mathbb{C}^{\text{op}} \rightarrow \mathcal{A}$ factors throught the inclusion $\mathcal{A}_{f.p} \overset{\iota}{\hookrightarrow} \mathcal{A}$ to give an equivalence $y|_{\mathcal{A}_{f.p}}: \mathbb{C}^{\text{op}} \rightarrow \mathcal{A}_{f.p}$.\\
    The composite of equivalences
    \vs 
    \begin{center}
        \begin{tikzcd}
            \mathbb{C}^{\text{op}} \arrow[rr, shift right=.8mm, bend right = 10mm, "\cong \epsilon_{\mathbb{C}}"]\arrow[r, "\iota^{0}"] & \mathcal{A}^{\text{op}} \arrow[r,"h_{\mathcal{A}}"] & \mathbf{LFC}\left(\mathcal{A}, \mathbf{Set}\right)
        \end{tikzcd}
    \end{center}
    \vs 
    is isomorphic to the component of the co-unit at $\mathbb{C}$, proving that the latter is also an equivalence. Notice here the notation $\iota^0:= \iota \circ y$.\\
\noindent \textbf{Theorem 1.1 (ii)} To prove the  second statement, observe that the inclusion
\begin{equation*}
    \mathcal{A}_{f.p} \simeq_{h_{\mathcal{A}}} \mathbf{LFC}\left(\mathcal{A}, \mathbf{Set}\right) \hookrightarrow \mathbf{CAT}\left(\mathcal{A}, \mathbf{Set}\right)
\end{equation*}
preserves finite limits. Therefore, $\iota: \mathcal{A}_{f.p} \hookrightarrow \mathcal{A}$ preserves finite colimits, since $h_{\mathcal{A}}$ is an equivalence of categories. Moreover, $\mathcal{A}_{f.p.}$ is essentially small, being the closure under colimits of a small set inside $\mathcal{A}_{f.p.}$. 
Then, by applying the functor F,
\begin{equation*}
        F\left(\mathcal{A}\right) = \textbf{LFC}\left(\mathcal{A}, \mathbf{Set}\right) \simeq \mathcal{A}_{f.p}^{\text{op}} \in \mathbf{Lex}
\end{equation*}
Observe that $h_{\mathcal{A}}$ induces an equivalence 
\begin{equation*}
    \mathbf{LEX}\left(\mathbf{LFC}\left(\mathcal{A}, \mathbf{Set}\right), \mathbf{Set}\right) \overset{h^{\ast}}{\rightarrow} \mathbf{LEX}\left(\mathcal{A}_{f.p.}^{\text{op}}, \mathbf{Set}\right)
\end{equation*}
that we can precompose with the component of the unit at $\mathcal{A}$
\begin{equation*}
    \eta_{\mathcal{A}}: \mathcal{A} \rightarrow \mathbf{LEX}\left(\mathbf{LFC}\left(\mathcal{A}, \mathbf{Set}\right), \mathbf{Set}\right)
\end{equation*}
to obtain 
\begin{align*}
    \tau:= h^{\ast} \circ \eta_{\mathcal{A}}: \mathcal{A} &\rightarrow \mathbf{LEX}\left(\mathcal{A}_{f.p.}^{\text{op}}, \mathbf{Set}\right)\\
A &\mapsto \mathcal{A}\left(\iota\left(-\right),A\right)
\end{align*}
To show: $\tau$ is an equivalence. If we prove that, it follows $\eta_{\mathcal{A}}$ is.   \\

For it, further machinery is needed.
\begin{definition}{[\cite{MP}, Definition 1.6]} A collection of objects $\mathcal{C}$ in a category $\mathcal{A}$ is conservative  iff, given a morphism $f:A \rightarrow B$ in $\mathcal{A}_1$ such that $\forall C \in \mathcal{C}_0$ the induced functor
    \begin{align*}
        f_{\ast}: \mathcal{A}\left(C,A\right) &\rightarrow \mathcal{A}\left(C,B\right)\\
        g &\mapsto fg 
    \end{align*}
    is a bijection, then f is an isomorphism. Therefore, if $\mathcal{A}$ is locally small, the collection $\mathcal{C}$ is conservative iff the Yoneda functor 
    \begin{align*}
        y: \mathcal{A} &\rightarrow \hat{\mathcal{C}}\\
        A &\mapsto \mathcal{A}\left(-, A\right)
    \end{align*}
    reflects isomorphisms.
\end{definition}
\begin{remark}
    The symbol $\hat{\mathcal{C}}$ denotes the categoy of presheaves over $\mathcal{C}$, i.e. $\hat{\mathcal{C}}:= \mathbf{Set}^{\mathcal{C}^{\text{op}}}$.
\end{remark}
\begin{lemma}{[\cite{MP}, Lemma 1.7]} Let $\mathcal{A}$ be a locally small category with filtered colimits, then 
    \begin{itemize}
        \item[(i)] Any set of generators in $\mathcal{A}$ is conservative;
        \item[(ii)] If $\mathcal{C}$ is a conservative set of f.p. objects in $\mathcal{A}$ which (as full subcategory) has finite colimits preserved by the inclusion $j: \mathcal{C} \hookrightarrow \mathcal{A}$, $j$ is dense. In other words, every object of $\mathcal{A}$ is a colimit of object of $\mathcal{C}$.
    \end{itemize}
\end{lemma}
\begin{corollary}{[\cite{MP}, Corollary 1.8]} If $\mathcal{A}$ is l.f.p., then the inclusion $\iota:\mathcal{A}_{f.p.} \hookrightarrow \mathcal{A}$ is dense.
\end{corollary}

\vs 

Back to the proof of the second part of Gabriel-Ulmer duality.\\

To show: $\tau$, given by
\begin{align*}
\tau: \mathcal{A} &\rightarrow \mathbf{LEX}\left(\mathcal{A}_{f.p.}^{\text{op}}, \mathbf{Set}\right)\\
A &\mapsto \mathcal{A}\left(\iota\left(-\right),A\right)
\end{align*}
is an equivalence.\\
By the previous lemma and [\cite{CWM}, page 247], $\iota$ is dense iff $A \mapsto \mathcal{A}\left(\iota \left(-\right), A \right)$ is fully faithful.\\
The last step consists in proving that $\tau$ is essentially surjective.\\  Let $M \in \mathbf{LEX}\left(\mathcal{A}_{f.p}^{\text{op}}, \mathbf{Set}\right)$, by density of the Yoneda embedding, M is representable as filtered colimit.
\begin{align*}
M &\cong \text{colim}_k \mathcal{A}_{f.p}^{\text{op}} \left(A_k, -\right) \cong \text{colim}\mathcal{A}\left(\iota \left(-\right), A_k \right) \cong  \\
&\cong \mathcal{A}\left(\iota \left(-\right), \text{colim}_k A_k \right)
\end{align*}
Thus: $M \cong \tau\left(A\right) $ for $A = \text{colim}_k A_k$. \ \ \ $\square$

\section{Application to Topos Theory}\label{appl}

\subsection{Elements of Topos Theory}\label{ele}

\ \ \ Let us assume the categories we are working with to be small, unless otherwise stated, and with finite left limits.
\begin{remark}
    Left limits are also called projective limits, especially in \cite{SGA}. The name left derives from the useful notation
    \begin{equation*}
        \lim_{\leftarrow I}
    \end{equation*}
    denoting such special limits.
\end{remark}

\vs  
In the first section we anticipated that the name left exact generalizes, in the context of topos theory, the analogous concept of homological algebra. We aim at providing evidence of that.

\begin{definition}{[\cite{MR}, 1, \S2, Defintion 2.4.1]} A site $\left(\mathbb{C}, \text{Cov(A)}\right)$ is a category $\mathbb{C}$ equipped with a Grothendieck topology on $\mathbb{C}$ given by a class of morphisms Cov(A) for each $A \in \mathbb{C}_0$. An element of Cov(A), the covering family of A, is an indexed set of morphisms over A, i.e.
    \begin{equation*}
    \left( A_i \overset{f_i}{\longrightarrow} A \right)_{i \in I}
    \end{equation*}
    The Grothendieck topology satisfies the following axioms:
    \begin{itemize}
        \item[(i)] Every isomorphism $f \in \mathbb{C}\left(A',A\right)$ gives a one-element covering family $\left\{f\right\} \in \text{Cov(A)}$;
        \item[(ii)] It is stable under pullbacks, i.e. in the following pullback diagram, 
        \vs
        \begin{center}
        \begin{tikzcd}
            A_i \times_B B \arrow[dr, phantom, very near start, "\lrcorner"]\arrow[d, swap, "g^{\ast}f_i"] \arrow[r] &A_i \arrow[d, "f_i"]\\
            B \arrow[r,swap, "g"] &A
        \end{tikzcd}
        \vs 
    \end{center} the pullback of the element $f_i$ of the covering family of A along $g$ is an element of the covering family of B. In other words, if $\left( f_i \right)_{i \in I} \in $ Cov(A) , then $g^{\ast}f_i \in$ Cov(B). 
\item[(iii)] It is closed under composition;
\item[(iv)] Monotonicity: If $\left(A_i \overset{f_i}{\longrightarrow} A \right)_{i \in I } \in $ Cov(A) and the sequentially composable family $\left(B_j \overset{g_j}{\longrightarrow}A \right)_{j \in J}$ is such that $\forall j \in J $ there exists $i \in I$ and a morphism $B_j \overset{s}{\longrightarrow} A_i$ making the following diagram commute 
\vs
\begin{center}
\begin{tikzcd}
    B_j \arrow[dr, swap,"g_j"]\arrow[r, "s"] &A_i \arrow[d,"f_i"]\\
    {} &A
\end{tikzcd}
\end{center}
\vs 
then $\left(g_j\right)_{j \in J} \in$ Cov(A).\\
This latter property allows the extracting of a subcover - in standard topological terms. 
\end{itemize}
\end{definition}
\begin{remark} The theory works also for locally small categories, by requiring the site $\mathbb{C}$ to have a topologically generating set $\mathcal{G}$: A set of objects $\left(A_i\right)_{i \in I}$ such that for each $A \in \mathbb{C}_0$ there exists a covering family $\left(A_i \overset{f_i}{\longrightarrow} A\right)_{i \in I}$.
\end{remark} 
    \begin{remark}There is an equivalent characterization of a Grothendieck topology in terms of sieves, subobjects in the category $\hat{\mathbb{C}}$ of presheaves over $\mathbb{C}$.
\end{remark}
\begin{example} The simplest example of a site is the one given by open sets of some topological space ordered by inclusion. Notice that a family of opens $\left(U_i \rightarrow U\right)_{i \in I} \in $ Cov(U) $\iff \ \bigcup_{i \in I} U_i = U$ with no additional local homeomorphism requirement. A site, therefore, generalizes the idea of covering, by getting rid of points.
\end{example}

\begin{definition}{[\cite{MR}, Definition 1.1.4]} Given two sites $\mathbb{C}, \mathbb{C}'$, we call a functor $\Psi: \mathbb{C} \rightarrow \mathbb{C}'$ a continuous functor iff it is left exact (i.e. preserves finite left limits) and it preserves coverings.
\end{definition}

\begin{remark}
    Such a functor is also called a $\mathbb{C}'$-model of $\mathbb{C}$. Indeed, in geometric logic (intuitionistic logic with $=, \land, \lor, \exists$ allowing for infinitary joins of formulae with finitely many variables), we have the following correspondence between geometric theories $\mathbb{T}$ and sites 
\begin{align*}
    \mathbb{T} &\leftrightarrow \left(\mathbb{C}, J\right) \ : \mathbb{C} \in \textbf{Lex}\\
    \text{Pt}\left(\mathcal{E}_{\mathbb{T}}\right) = \text{Sh}\left(\mathbb{C}, J \right) \simeq \text{Mod}\left(\mathbb{T}\right) &\hookrightarrow \text{LEX}\left(\mathbb{C}, \mathbf{Set}\right) 
\end{align*}
where $\mathcal{E}_{\mathbb{T}}$ is the classifying topos of the theory $\mathbb{T}$.
\end{remark}
\begin{definition} A presheaf over a category $\mathbb{C}$ is a contravariant functor into $\mathbf{Set}$.
\end{definition}

\begin{definition}{[\cite{MR}, 1, \S1, Definition 1.1.6]}
Given a presheaf $F: \mathbb{C} \rightarrow \mathbf{Set}$, a compatible family of morphisms from a covering $\left(A_i \overset{f_i}{\rightarrow}A \right)_{i \in I}$ to the presheaf $F$ is a family of morphisms
\begin{equation*}
    \left(A_i \overset{\xi_i}{\rightarrow}F\right)_{i \in I} 
\end{equation*}
such that the external square in the diagram below commutes
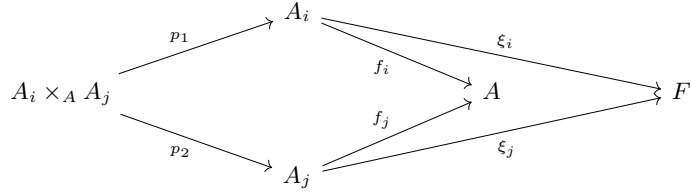
\begin{center}
    \begin{tikzcd}[column sep=2cm]
        {} &A_i \arrow[dr,swap, "f_i"] \arrow[drr,"\xi_i"]&{} &{}\\
        A_i \times_A A_j \arrow[ur,"p_1"]\arrow[dr, swap,"p_2"] &{} &A &F\\
        {} &A_j \arrow[ur, "f_j"] \arrow[urr,swap,"\xi_j"] &{} &{}
    \end{tikzcd}
\end{center}
\vs 
i.e.
\begin{equation*}
    \xi_jp_2 = \xi_i p_1
\end{equation*}
for all indices $i,j \in I$.
\end{definition} 
\begin{definition}{[\cite{MR}, ibid.]} A presheaf $F$ is a sheaf for $\mathbb{C}$ iff for $A \in \mathbb{C}_0,$ and $\left(\xi\right)_{i \in I}$ a compatible family for the covering $\left(f_i\right)_{i \in I}$, there exists a unique morphism $\xi:A \rightarrow F$ such that $ \xi f_i = f, i \in I$.
    \begin{center}
        \begin{tikzcd}[column sep=2cm]
            {} &A_i \arrow[dr,swap, "f_i"] \arrow[drr,"\xi_i"]&{} &{}\\
            A_i \times_A A_j \arrow[ur,"p_1"]\arrow[dr, swap,"p_2"] &{} &A \arrow[r, near start, "! \xi"] &F\\
            {} &A_j \arrow[ur, "f_j"] \arrow[urr,swap,"\xi_j"] &{} &{}
        \end{tikzcd}
    \end{center}
\end{definition}
\begin{definition} {[\cite{MR}, ibid.]} A presheaf $F$ over a site $\mathbb{C}$ is called separated presheaf iff there exists at most one $\xi$ as defined above.
\end{definition}
Finally, we are able to see where exactness comes into play in defining a sheaf, by analyzing the case of $\mathbf{Set}$.
\begin{example} In $\mathbf{Set}$, for a presheaf $F$ and a family $\left(f_i\right)_{i \in I}$, F is a sheaf if the following diagram is left exact
    \vs 
    \begin{center}
        \begin{tikzcd}
            \hat{\mathbb{C}}\left(A, F\right) \arrow[r,"u"] & \prod_{i \in I} \hat{\mathbb{C}}\left(A_i, F \right) \arrow[r, shift right, swap,"v_2"] \arrow[r, shift left, "v_1"] &\prod_{i,j \in I} \hat{\mathbb{C}}\left(A_i \times_A A_j, F\right)
        \end{tikzcd}
    \end{center}
    \vs
    In fact, $F$ being a sheaf amounts to 
    \begin{equation*}
        u\left(\xi\right)= \langle \xi f_i|i \in I\rangle
    \end{equation*}
    being the equalizer of $v_1, v_2$, i.e. a monic arrow, and $\left(f_i\right)_{i \in I}$ being a covering family.\\
    Further, the maps $v_1, v_2$ are defined as follows:
    \begin{align*}
        v_1\left(\langle h_i:i \in I\rangle \right)&= \langle h_i \circ p_1^{i,j}: i,j \in I \rangle \\
        v_2 \left(\langle h_i: i \in I \rangle \right) &= \langle h_i \circ p_2^{i,j}: i,j \in I \rangle 
    \end{align*}
    \end{example}

    The category of sheaves over $\mathbb{C}$, denoted by $\tilde{\mathbb{C}}$, is a full subcategory of $\hat{\mathbb{C}}$. Furthermore, there exists a canonical adjunction
    \vs 
    \begin{center}
        \begin{tikzcd}
        \tilde{\mathbb{C}} \arrow[r,hook, shift left,"\iota"] &\hat{\mathbb{C}} \arrow[l,shift left, "a"]
        \end{tikzcd}
    \end{center}
    where $\iota \vdash a$ and $a$ preserves filtered colimits. For $F \in \hat{\mathbb{C}}, a\left(F\right)$ is the so-called associated sheaf (to F).

\begin{definition}{[\cite{MR}, 1, \S3, Definition 1.3.1]} A Grothendieck topos $\mathcal{E}$ is a category equivalent to the one of sheaves over a small site $\mathbb{C}$.
    \begin{equation*}
        \mathcal{E} \simeq \tilde{\mathbb{C}}
    \end{equation*}
    \end{definition}

    \subsection{Global Sections and Points}

    \ \ \ In geometry, dealing with (pre)sheaves on a space, we would like to be able to reconstruct global data from a set of local ones.
    For it, we study the interplay between stalks of the presheaf and global sections.\\

    If $\mathbb{C}$ is a small, finitely complete category and $\mathcal{E}$ is a Grothendieck topos, a functors
    \begin{equation*}
        F  \in \mathbf{LEX}\left(\mathbb{C}, \mathcal{E}\right)
    \end{equation*}
    provides a sheaf of models of $\mathbb{C}$ over the generalized space $\mathcal{E}$.\\

   The sheaf functor $\epsilon: \mathbb{C} \rightarrow \tilde{\mathbb{C}}$ is left exact and continuous. For all $\mathcal{E}$ Grothendieck topoi and $\mathcal{E}$-model $F$ of $\mathbb{C}$ there exists a n $\mathcal{E}$-model $\tilde{F}$ of $\tilde{\mathbb{C}}$ such that the diagram commutes
   \vs
   \begin{center}
    \begin{tikzcd}
        \mathbb{C} \arrow[dr, swap, "F"] \arrow[r, "\epsilon"] &\tilde{\mathbb{C}} \arrow[d, "\tilde{F}"]\\
        {} &\mathcal{E}
    \end{tikzcd}
\end{center}
\vs
\begin{definition}{[\cite{SGA},4,\S3,Definition 3.1]} A morphism of topoi, also called geometric morphism, is a triple of functors 
    \begin{equation*}
        u=\left(u_{\ast}, u^{\ast}, \phi\right)
    \end{equation*}
    where $u_{\ast} \vdash u^{\ast}$ by the adjunction
    \begin{equation*}
        \phi: \mathcal{E} \left( u^{\ast}X', Y \right) \overset{\sim}{\rightarrow} \mathcal{E}'\left(X', u_{\ast}Y\right)
    \end{equation*}
    for $X' \in \mathcal{E}_0', X \in \mathcal{E}_0$.
\end{definition} 
\begin{definition} We denote by $\mathbf{TOP}$ the 2-category of Grothendieck topoi and geometric morphisms between them.
\end{definition}

Having introduced the notion of a sheaf of models, define the global sections of $F$ by postcomposition with 
\begin{equation*}
    \Gamma := \mathcal{E}\left(1, - \right) : \mathcal{E} \rightarrow \mathbf{Set}
\end{equation*}
Namely,
\vs 
\begin{center}
    \begin{tikzcd}
        \mathbb{C} \arrow[rr, shift right=.8mm, bend right = 10mm, swap,"\Gamma F"] \arrow[r,"F"] &\mathcal{E} \arrow[r, "\Gamma"] & \mathbf{Set}
    \end{tikzcd}
\end{center}
\vs
\begin{equation*}
    \Gamma F \in \mathbf{LEX}\left(\mathbb{C}, \mathbf{Set}\right)
\end{equation*}
\vs
Dually, given a point of $\mathcal{E}$, i.e. a geometric morphism 
\begin{equation*}
    p \in \mathbf{TOP}\left(\mathbf{Set}, \mathcal{E}\right) = \text{Pt}\left(\mathcal{E}\right)
\end{equation*}
the stalk of $F$ at $p$, denoted by $F_p$, is also a left exact functor
\vs 
\begin{center}
    \begin{tikzcd}
        \mathbb{C} \arrow[rr, shift right=.8mm, bend right = 10mm, swap,"F_p"] \arrow[r,"F"] &\mathcal{E} \arrow[r, "p^{\ast}"] & \mathbf{Set}
    \end{tikzcd}
\end{center}
\vs where $p^{\ast}$ is the inverse image part of $p$.\\

Taking $\mathbf{Set}$-models of the site, $F$ is a sheaf representation of $\Gamma F$ via the stalk of models $F_p$. Can we do the whole business only in $\mathbf{LEX}\left(\mathbb{C}, \mathbf{Set}\right)?$\\
\begin{example}
For a presheaf topos over a small site, $\Gamma : \mathcal{E} \rightarrow \mathbf{Set}$ is a limit in $\mathbb{C}^{\text{op}}$ and for each object of the site 
\begin{equation*}
    c \in \mathbb{C} \rightarrow p_c \in \text{Pt}\left(\mathcal{E}\right)= \text{Pt}\left( \mathbf{Set}^{\mathbb{C}^{\text{op}}}\right)
\end{equation*}
and its inverse part is given by evaluation 
\begin{equation*}
    p_c^{\ast}\left(F\right) = F\left(c\right)
\end{equation*}
where $F \in \mathbf{LEX}\left(\mathbb{C}, \mathbf{Set}^{\mathbb{C}^\text{op}}\right)$ is the sheaf of models of $\mathbb{C}$ over the presheaf topos.\\
Moreover, a change of object in the site $\gamma: c \rightarrow c'$ gives rise to a natural transformation 
\begin{equation*}
    p_{c'}^{\ast} \rightarrow p_c
\end{equation*}
(in components) that allows to have a diagram whose limit gives $\Gamma$ 
\begin{equation*}
    \Gamma \cong \lim_{c \in \mathbb{C}^{\text{op}}} p_c^{\ast} 
\end{equation*}
Thus: taken $F \in \mathbf{LEX}\left(\mathbb{C}, \mathbf{Set}^{\mathbb{C}^\text{op}} \right)$, the global sections of F is obtained as limit of stalks of F and this limit is in $\mathbf{LEX}\left(\mathbb{C}, \mathbf{Set}\right)$.
\end{example}

\subsection{Gabriel-Ulmer Duality for Topoi}

\ \ \ For general topoi it is not always the case that from stalks it is possible to reconstruct the global sections via a limit procedure. However, when a topos has enough points, such property still holds.  
\begin{definition}
    A Grothendieck topos $\mathcal{E}$ has enough points iff for $f:X \rightarrow Y $ in $\mathcal{E}$ and $p^{\ast}\left(f\right)$ being a bijection $\forall p \in \text{Pt}\left(\mathcal{E}\right)$, $f$ is an isomorphism.
\end{definition}
\begin{definition} For a collection of objects $\mathcal{K}$ in $\mathcal{A} \in \mathbf{LFP}_0$, a locally finitely presentable category, we denote by $\left[\mathcal{K}\right]$ the least full subcategory of $\mathcal{A}$ containing $\mathcal{K}$ closed under limits and filtered colimits.
\end{definition}
\begin{remark}
    The category just defined is locally finitely presentable.
\end{remark}
\begin{theorem}{[\cite{MP}, Theorem 3.3]} If $\mathbb{C}$ is a small category with finite limits, $\mathcal{E}$ has enough points and $F \in \mathbf{LEX}\left(\mathbb{C}, \mathbf{Set}\right)$. For $L \in \mathbf{LEX}\left(\mathcal{E}, \mathbf{Set}\right)$
    \begin{equation*}
        LF \in \left[ \left\{F_p | p \in \text{Pt}\left(\mathcal{E}\right)\right\}\right] \hookrightarrow \mathbf{LEX} \left( \mathbb{C}, \mathbf{Set}\right) 
    \end{equation*}
    In particular, taking $L = \Gamma$, the global sections of $F$ is the closure of the stalks of $F$ under limits and filtered colimits.
\end{theorem}
\begin{remark} Notice that by the considerations made at the end of the previous section $F \in \mathbf{LEX}\left(\mathbb{C}, \mathbf{Set}\right)$.
\end{remark}

Therefore, Gabriel-Ulmer duality relates the site $\mathbb{C} \in \mathbf{Lex}$ to the locally finitely presentable $\mathbf{LFC}$-closed category of stalks of F and allows, when $\mathcal{E}$ has enough points, to go back and forth between local and global data.

\section{The geometry of clans: Future directions}\label{future}

\ \ \ In [\cite{JF}, Definition 2.1] a clan is defined as a small category $\mathbb{C}$ with a terminal object 1 equipped with a class of maps, called display maps such that 
\begin{itemize}
\item[(i)] Pullbacks of display maps are display maps;
\item[(ii)] The class of display maps is closed under composition;
\item[(iii)] Isomorphisms are display maps;
\item[(iv)] Terminal projections are display maps.
\end{itemize}

It is apparent the analogy with the axioms defining a Grothendieck topology. In fact, suppose that given two objects $A,B \in \mathbb{C}_0$, there exists an arrow $A \overset{r}{\rightarrow} B$ such that the following diagram commutes
\vs 
\begin{center}
\begin{tikzcd}
A \arrow[dr, swap, "!f"] \arrow[r, "r"] &B \arrow[d, "!g"]\\
{} &1
\end{tikzcd}
\end{center}
\vs
then the composite $g r$ is also a display map. This condition replaces the one given in the context of Grothendieck topologies for monotonicity.\\

Moreover, in \cite{JF} Gabriel-Ulmer duality is extended to a duality between Cauchy-complete clans and clan algebraic categories. [\cite{JF}, (1.1)]:
\begin{equation}
\mathbf{Clan}_{\text{cc}} \simeq \mathbf{ClanAlg}^{\text{op}}
\end{equation}

Future work will include the investigation of the relation between clans and topos theory, as clans generalize the notion of site and display maps the one of coverings, mutatis mutandis.



\newpage 
\begin{thebibliography}{10}
        \bibitem[ARV]{ARV} Adámek, J., Rosický, J., Vitale, E., Algebraic Theories, CUP, 2011;
	\bibitem[MP]{MP} Makkai, M., Pitts, A.,M., Some Results on Locally Finitely Presentable Categories, Transactions of the AMS, Volume 299, Num.2, February 1987;
	\bibitem[CWM]{CWM} Mac Lane, S., Categories for the Working Mathematician, Second Edition, Graduate Texts in Mathematics, vol. 5, Springer-Verlag, New York, Inc., 1978;
    \bibitem[MR]{MR} Makkai, M., Reyes, G., E., First Order Categorical Logic, Lecture Notes in Mathematics, vol. 611, Springer-Verlag, Berlin, Heidelberg, New York, 1977; 
    \bibitem[AR]{AR} Adámek, J., Rosický, J., Locally Presentable and Accessible Categories, London Mathematical Society Lecture Note Series, vol. 189, Cambridge University Press, 1994;
    \bibitem[SGA4.I]{SGA} Théorie des Topos et Cohomologie Etale des Schémas, Tome 1, Théorie des Topos, Eds., M. Artin, A. Grothendieck and J.L. Verdier, Lecture Notes in Math., vol. 269, Springer-Verlag, Berlin and New York, 1972;
    \bibitem[GU]{GU} Gabriel, P., Ulmer, F., Lokal praesentierbare Kategorien, Lecture Notes in Math., vol. 221, Spinger-Verlag, Berlin and New York, 1971;
    \bibitem[JF]{JF} Frey, J., Duality for Clans: An extension of Gabriel-Ulmer Duality, arXiv e-prints, 2023, arXiv:2308.11967;
    \bibitem[DLRG]{DLRG} Di Liberti, I., Ramos González, J., Gabriel–Ulmer Duality for Topoi and its Relation with Site Presentations. Applied Categorical Structures 28, 935–962, Springer, 2020, https://doi.org/10.1007/s10485-020-09605-x;
    \bibitem[KE]{KE} Category Seminar, Proceedings Sydney Category Theory Seminar 1972 /1973, edited by Gregory M. Kelly, Lecture Notes in Mathematics, 420, Review of the elements of 2-categories by Kelly, G., M., and Street, R., Springer Verlag Berlin, Heidelberg, 1974, https://doi.org/10.1007/BFb0063096;
    \bibitem[PV]{PV} Palmgren, E., Vickers, S., J., Partial Horn logic and cartesian categories, Annals of Pure and Applied Logic, 145, 314-355, 2007.


\end{thebibliography}

\end{document}